\documentclass[11pt]{article}
\usepackage{amscd}
\usepackage{amsmath}
\usepackage{latexsym}
\usepackage{amsfonts}
\usepackage{amssymb}
\usepackage{amsthm}
\usepackage{graphicx}

\newtheorem{theorem}{Theorem}[section]
\newtheorem{lemma}{Lemma}[section]

\newtheorem{proposition}{Proposition}[section]

\newtheorem{remark}{Remark}[section]

\makeatletter
\newcommand{\Extend}[5]{\ext@arrow0099{\arrowfill@#1#2#3}{#4}{#5}}
\makeatother

\begin{document}

 \title{ Global well-posedness and scattering for the defocusing $H^{\frac12}$-subcritical Hartree equation in $\mathbb{R}^d$}
 \author{{Changxing Miao,\ \ Guixiang Xu,\ \ and \ Lifeng Zhao }\\
         {\small Institute of Applied Physics and Computational Mathematics}\\
         {\small P. O. Box 8009,\ Beijing,\ China,\ 100088}\\
         {\small (miao\_changxing@iapcm.ac.cn, \ xu\_guixiang@iapcm.ac.cn, zhao\_lifeng@iapcm.ac.cn ) }\\
         \date{}
        }
\maketitle

\begin{abstract}\noindent
We prove the global well-posedness and scattering for the defocusing
$H^{\frac12}$-subcritical (that is, $2<\gamma<3$) Hartree equation
with low regularity data in $\mathbb{R}^d$, $d\geq 3$. Precisely, we
show that a unique and global solution exists for initial data in
the Sobolev space $H^s\big(\mathbb{R}^d\big)$ with
$s>4(\gamma-2)/(3\gamma-4)$, which also scatters in both time
directions. This improves the result in \cite{ChHKY}, where the
global well-posedness was established for any
$s>\max\big(1/2,4(\gamma-2)/(3\gamma-4)\big)$. The new ingredients
in our proof are that we make use of an interaction Morawetz
estimate for the smoothed out solution $Iu$, instead of an
interaction Morawetz estimate for the solution $u$, and that we make
careful analysis of the monotonicity property of the multiplier
$m(\xi)\cdot \langle \xi\rangle^p$. As a byproduct of our proof, we
obtain that the $H^s$ norm of the solution obeys the uniform-in-time
bounds.
\end{abstract}

 \begin{center}
 \begin{minipage}{120mm}
 \end{minipage}
 \end{center}


\section{Introduction}
 \setcounter{section}{1}\setcounter{equation}{0}
 In this paper, we study the global well-posedness of the following
 initial value problem (IVP) for the
defocusing $H^{\frac12}$-subcritical (that is, $2<\gamma<3$) Hartree
equation.
\begin{equation} \label{equ1}
\left\{ \aligned
    iu_t +  \Delta u  & = \big( |x|^{-\gamma}* |u|^2 \big) u, \quad d\geq 3,\\
     u(0) & =u_0(x)\in   H^s(\mathbb{R}^d),
\endaligned
\right.
\end{equation}
where $H^s$ denotes the usual inhomogeneous Sobolev space of order
$s$.

We adopt the following standard notion of local well-posedness, that
is, we say that the IVP (\ref{equ1}) is locally well-posed in $H^s$
if for any $u_0 \in H^s$, there exists a positive time
$T=T(\big\|u_0\big\|_{s})$ depending only on the norm of the initial
data, such that a solution to the IVP exists on the time interval
$[0, T]$, is unique in a certain Banach space  of functional $X
\subset C\big([0, T], H^s \big)$, and the solution map from $H^s_x$
to $C\big([0, T], H^s \big)$ depends continuously. If $T$ can be
taken arbitrarily large, we say that the IVP (\ref{equ1}) is
globally well-posed.

Local well-posedness for the IVP (\ref{equ1}) in $H^s$ for any
$s>\frac{\gamma}{2}-1$ was established in \cite{MiXZ06}. A local
solution also exists for $H^{\frac{\gamma}{2}-1}$ initial data, but
the time of existence depends not only on the
$H^{\frac{\gamma}{2}-1}$ norm of $u_0$, but also on the profile of
$u_0$. For more details on local well-posedness see \cite{MiXZ06}.

$L^2$ solutions of (\ref{equ1}) enjoy mass conservation
\begin{equation*}
\aligned \big\|u(t,\cdot)\big\|_{L^2\big(
\mathbb{R}^d\big)}=\big\|u_0(\cdot)\big\|_{L^2\big(
\mathbb{R}^d\big)}.
\endaligned
\end{equation*}
Moreover, $H^1$ solutions enjoy energy conservation
\begin{equation*}
\aligned
 E(u)(t)=\frac12 & \big\|\nabla
u(t)\big\|^2_{L^2\big( \mathbb{R}^d\big)}+\frac{1}{4}
\iint_{\mathbb{R}^d\times \mathbb{R}^d} \frac{1}{|x-y|^{\gamma}}
|u(t,x)|^2 |u(t,y)|^2\ dxdy=E(u)(0),
\endaligned
\end{equation*}
which together with mass conservation and the local theory
immediately yields global well-posedness for (\ref{equ1}) with
initial data in $H^1$. A large amount of works have been devoted to
global well-posedness and scattering for the Hartree equation, see
\cite{GiO93}-\cite{GiV01}, \cite{HaT87}, \cite{LiMZ08},
\cite{Mi97}-\cite{NaO92}.

\begin{figure}[ht]\label{figure}
\centering
\includegraphics[width=0.6\textwidth, angle=0]{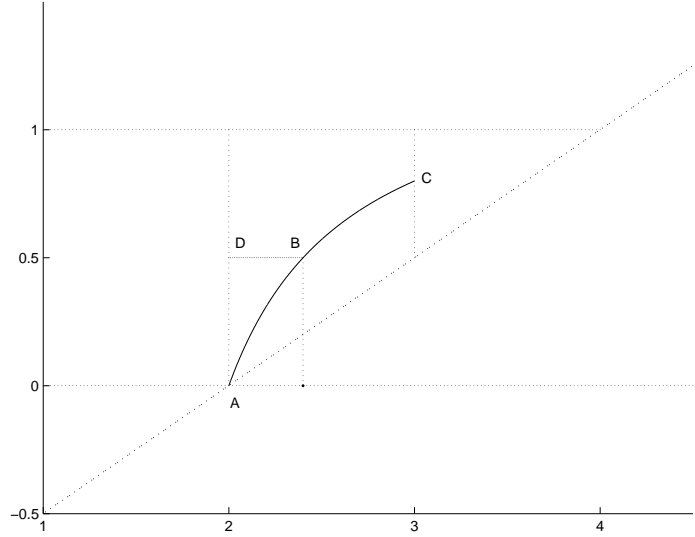}
\caption[]{The curve ``ABC'' is descripted by ``$s=\frac{4(\gamma-2)}{3\gamma-4}$''}
\end{figure}

Existence of global solutions in $\mathbb{R}^3$ to (\ref{equ1})
corresponding to initial data below the energy threshold was
recently obtained in \cite{ChHKY} by using the method of ``almost
conservation laws'' or ``I-method'' (for a detailed description of
this method, see \cite{Taobook1} or section 3 below) and the
interaction Morawetz estimate for the solution $u$, where global
well-posedness was obtained in $H^s(\mathbb{R}^3)$ with
$s>\max\big(1/2,4(\gamma-2)/(3\gamma-4)\big)$. Since authors in
\cite{ChHKY} used the interaction Morawetz estimate,  which involves
$\dot{H}^{1/2}$ norm of the solution, the restriction condition $s
\geq \frac{1}{2}$ is prerequisite. In order to resolve IVP
(\ref{equ1}) in $H^s$, $s<\frac12$ by still using the interaction
Morawetz estimate, we need return to the interaction Morawetz
estimate for the smoothed out version $Iu$ of the solution, which is
initially used in \cite{CGT07}, whereafter in \cite{DPST07}.

In this paper, we consider the case $d\geq 3$ and we prove the
following result:

\begin{theorem}\label{theorem}
Let $2<\gamma<3 \leq d$, the initial value problem (\ref{equ1}) is
globally well-posedness in $H^s(\mathbb{R}^d)$ for any
$s>\frac{4(\gamma-2)}{3\gamma-4}$. Moreover the solution satisfies
\begin{equation*}
\aligned \sup_{t\in[0, \infty)}\big\|
u(t)\big\|_{H^s\big(\mathbb{R}^d\big)} \leq
C\big(\big\|u_0\big\|_{H^s}\big),
\endaligned
\end{equation*}
and there is scattering for these solutions, that is, the wave
operators exist and there is asymptotic completeness on all of
$H^s(\mathbb{R}^d)$.
\end{theorem}
\begin{remark}
As for the case $3\leq \gamma < 4 \leq d$, local well-posedness for
the IVP (\ref{equ1}) in $H^s$ holds for any $s>\frac{\gamma}{2}-1$.
Note that in this case, we have
\begin{equation*}
\aligned \frac{\gamma}{2}-1 \geq \frac12,
\endaligned
\end{equation*}
which satisfies the need of the regularity of the interaction
Morawetz estimate. Hence we only combine ``I-method'' with the
interaction Morawetz estimate for the solution to obtain the low
regularity of the IVP (\ref{equ1}), just as in \cite{CKSTT04} .
\end{remark}

For the case $d=3$, Theorem \ref{theorem} improves the result
$s>\max\big(1/2,4(\gamma-2)/(3\gamma-4)\big)$ in \cite{ChHKY} (see
Figure \ref{figure}), where the authors used ``I-method'' and the
interaction Morawetz estimate for the solution just as in
\cite{CKSTT04}. In general, in order to prove the almost
conservation law, one doesn't need to use the monotonicity property
of the multiplier $m(\xi)\cdot \langle \xi\rangle^p$. In the present
paper, we prove Theorem \ref{theorem} by combining I-method with an
interaction Morawetz estimate for the smoothed out version $Iu$ of
the solution. Such a Morawetz estimate for an almost solution is the
main novelty of this paper, which can lower the need on the
regularity of the initial data.

Last, we organize this paper as following: In Section 2, we
introduce some notation and state some important propositions that
we will used throughout this paper. In Section 3, we review the
I-method, prove the local well-posedness theory for $Iu$ and obtain
an upper bound on the increment of the modified energy. In Section
4, we prove the ``almost interaction Morawetz estimate'' for the
smoothed out version $Iu$ of the solution. Finally in Section 5, we
give the details of the proof of the global well-posedness stated in
Theorem \ref{theorem}.

\section{Notation and preliminaries}
 \setcounter{section}{2}
\setcounter{equation}{0}

\subsection{Notation} In what follows, we use $A\lesssim B$ to denote
an estimate of the form $A \leq CB$ for some constant $C$. If
$A\lesssim B$ and $B\lesssim A$, we say that $A\thickapprox B$. We
write $A\ll B$ to denote an estimate of the form $A \leq cB$ for
some small constant $c>0$. In addition $\langle a \rangle:=1+|a|$
and $a\pm:=a\pm \epsilon$ with $0<\epsilon \ll 1$. The reader also
has to be alert that we sometimes do not explicitly write down
constants that depend on the $L^2$ norm of the solution. This is
justified by the conservation of the $L^2$ norm.

\subsection{Definition of spaces}
We use $L^r_x(\mathbb{R}^d)$ to denote the Lebesgue space of
functions $f:\mathbb{R}^d\rightarrow \mathbb{C}$ whose norm
\begin{equation*}
\aligned \big\|f\big\|_{L^r_x}:=\Big(\int_{\mathbb{R}^d} \big| f(x)
\big|^r dx\Big)^{\frac1r}
\endaligned
\end{equation*}
is finite, with the usual modification in the case $r=\infty$. We
also use the space-time Lebesgue spaces $L^q_tL^r_x$ which are
equipped with the norm
\begin{equation*}
\aligned \big\|u\big\|_{L^q_tL^r_x}:=\Big( \int_{J}
\big\|u(t,x)\big\|^{q}_{L^r_x} dt \Big)^{\frac1q}
\endaligned
\end{equation*}
for any space-time slab $J\times \mathbb{R}$, with the usual
modification when either $q$ or $r$ are infinity. When $q=r$, we
abbreviate $L^q_tL^r_x$ by $L^q_{t,x}$.

As usual, we define the Fourier transform of $f(x)\in L^1_x$ by
\begin{equation*}
\aligned
\widehat{f}(\xi)=\big(2\pi\big)^{-\frac{d}{2}}\int_{\mathbb{R}^d}
e^{-ix\xi} f(x)dx.
\endaligned
\end{equation*}
We define the fractional differentiation operator
$|\nabla_x|^{\alpha}$ for any real $\alpha$ by
\begin{equation*}
\aligned \widehat{|\nabla|^{\alpha}u}(\xi):=|\xi|^\alpha
\widehat{u}(\xi),
\endaligned
\end{equation*}
and analogously
\begin{equation*}
\aligned \widehat{\langle \nabla\rangle^{\alpha}u}(\xi):=\langle \xi
\rangle^\alpha \widehat{u}(\xi).
\endaligned
\end{equation*}
The inhomogeneous Sobolev space $H^s(\mathbb{R}^d)$ is given via
\begin{equation*}
\aligned \big\|u\big\|_{H^s}:=\big\|\langle \nabla\rangle^s
u\big\|_{L^2\big(\mathbb{R}^d\big)},
\endaligned
\end{equation*}
while the homogeneous Sobolev space $\dot{H}^s(\mathbb{R}^d)$ is
given via
\begin{equation*}
\aligned \big\|u\big\|_{\dot{H}^s}:=\big\|| \nabla|^s
u\big\|_{L^2\big(\mathbb{R}^d\big)}.
\endaligned
\end{equation*}

Let $S(t)$ denote the solution operator to the linear
Schr\"{o}dinger equation
\begin{equation*}
\aligned iu_t+\Delta u=0, x\in \mathbb{R}^d.
\endaligned
\end{equation*}
We denote by $X^{s,b}(\mathbb{R}\times \mathbb{R}^d)$ the completion
of $\mathcal{S}\big(\mathbb{R}\times \mathbb{R}^d\big)$ with respect
to the following norm
\begin{equation*}
\aligned
\big\|u\big\|_{X^{s,b}}=\big\|S(-t)u\big\|_{H^s_xH^b_t}=\big\|\langle
\tau+|\xi|^2 \rangle^{b} \langle \xi \rangle^{s}
\widetilde{u}(\tau,\xi)\big\|_{L^2_{\tau}L^2_{\xi}\big(\mathbb{R}\times
\mathbb{R}^d\big)},
\endaligned
\end{equation*}
where $\widetilde{u}$ is the space-time Fourier transform
\begin{equation*}
\aligned \widetilde{u}(\tau,\xi)=\big(2
\pi\big)^{-\frac{d+1}{2}}\iint_{\mathbb{R}\times \mathbb{R}^d}
e^{-i(x\cdot \xi + t\tau)}u(t,x)dtdx.
\endaligned
\end{equation*}
Furthermore for a given time interval $J$, we define
\begin{equation*}
\aligned \big\|u\big\|_{X^{s,b}(J)}=\inf \Big\{
\big\|v\big\|_{X^{s,b}};\ v=u \ \text{on}\ J\Big\}.
\endaligned
\end{equation*}

\subsection{Some known estimates}
Now we recall a few known estimates that we shall need. First we
state the following Strichartz estimate \cite{Ca03}, \cite{KeT98}.
Let $d\geq 3$, we recall that a pair of exponents $(q,r)$ is called
admissible if
\begin{equation*}
\aligned \frac{2}{q}=d(\frac12-\frac1r),\quad 2\leq q, r \leq
\infty.
\endaligned
\end{equation*}

\begin{proposition}\label{stri}
Let $d\geq 3$,  $(q,r)$ and $(\widetilde{q},\widetilde{r})$ be any
two admissible pairs. Suppose that $u$ is a solution to
\begin{equation*}
\aligned iu_t + \Delta u & = F(t,x), \ t\in J, x\in \mathbb{R}^d \\
u(0) & = u_0(x).
\endaligned
\end{equation*}
Then we have the estimate
\begin{equation*}
\aligned \big\|u\big\|_{L^q_tL^r_x\big(J\times \mathbb{R}^d\big)}
\lesssim \big\|u_0\big\|_{L^2\big( \mathbb{R}^d \big)} +
\big\|F\big\|_{L^{\widetilde{q}'}_tL^{\widetilde{r}'}_x\big(J\times
\mathbb{R}^d \big)},
\endaligned
\end{equation*}
where the prime exponents denote H\"{o}lder dual exponents.
\end{proposition}

Let us say that a function $u$ has spatial frequency $N$ if its
Fourier transform is supported on the annulus $\big\{\langle
\xi\rangle \thickapprox N\big\}$. From Strichartz estimate
$\big\|u\big\|_{L^{q}_tL^r_x} \lesssim \big\|u\big\|_{X^{0,
\frac12+}}$ for admissible $(q,r)$ and Sobolev embedding theorem, we
have
\begin{proposition}\label{linear}
For $r<\infty, \ 0 \leq \frac{2}{q} \leq \min\big( \delta(r), 1
\big)$, we have
\begin{equation*} \aligned
\big\|u\big\|_{L^q_tL^r_x} \lesssim
\big\|u\big\|_{X^{\delta(r)-\frac2q , \frac12+}}.
\endaligned
\end{equation*}
While for $2\leq q \leq \infty, r=\infty$, we have
\begin{equation*} \aligned
\big\|u\big\|_{L^q_tL^{\infty}_x} \lesssim
\big\|u\big\|_{X^{\frac{d}{2}-\frac2q+ , \frac12+}}.
\endaligned
\end{equation*}
\end{proposition}

%

\section{the I-method and the modified local well-posedness}
\setcounter{section}{3} \setcounter{equation}{0}

\subsection{the I-operator and the hierarchy of energies}
Let us define the operator $I$. For $s<1$ and a parameter $N\gg 1$,
let $m(\xi)$ be the following smooth monotone multiplier:
\begin{equation*}
\aligned m(\xi):=\left\{\begin{array}{ll} 1, & \text{if}\ \ |\xi|<N,
\\
\big( \frac{N}{|\xi|}\big)^{1-s}, & \text{if}\ \ |\xi|>2N.
\end{array}\right.
\endaligned
\end{equation*}
We define the multiplier operator $I:H^s \longrightarrow H^1$ by
\begin{equation*}
\aligned \widehat{Iu}(\xi) = m(\xi)\widehat{u}(\xi).
\endaligned
\end{equation*}
The operator $I$ is smoothing of order $1-s$ and we have that
\begin{equation*}
\aligned
 \big\|u\big\|_{H^{s_0}} \lesssim &
\big\|Iu\big\|_{H^{s_0+1-s}}\lesssim
N^{1-s}\big\|u\big\|_{H^{s_0}},\\
\big\|u\big\|_{X^{s_0, b_0}} \lesssim &\big\|Iu\big\|_{X^{s_0+1-s,
b_0}}\lesssim N^{1-s}\big\|u\big\|_{X^{s_0, b_0}}
\endaligned
\end{equation*}
for any $s_0, b_0\in \mathbb{R}$.

We set
\begin{equation}\label{I1}
\aligned \widetilde{E}(u)=E(Iu),
\endaligned
\end{equation}
where
\begin{equation*}
\aligned E(u)(t)=\frac12 & \big\|\nabla
u(t)\big\|^2_{L^2}+\frac{1}{4} \int \int \frac{1}{|x-y|^{\gamma}}
|u(t,x)|^2 |u(t,y)|^2\ dxdy.
\endaligned
\end{equation*}
We call $\widetilde{E}(u)$ the modified energy. Since we will focus
on the analysis of the modified energy, we collect some facts
concerning the calculus of multilinear forms used to define the
modified energy.

If $k\geq 2$ is an even integer, we define a spatial multiplier of
order $k$ to be the function $M_k(\xi_1, \xi_2,\cdots, \xi_k)$ on
\begin{equation*}
\aligned \Gamma_k=\Big\{(\xi_1, \xi_2,\cdots, \xi_k)\in
\big(\mathbb{R}^d\big)^k: \sum^k_{j=1}\xi_j=0 \Big\},
\endaligned
\end{equation*}
which we endow with the standard measure
$\delta(\xi_1+\xi_2+\cdots+\xi_k)$. If $M_k$ is a multiplier of
order $k$, $1\leq j\leq k$ is an index and $l\geq 1$ is an even
integer, the elongation $X^l_j(M_k)$ of $M_k$ is defined to be the
multiplier of order $k+l$ given by
\begin{equation*}
\aligned X^l_j(M_k)(\xi_1,\xi_2,\cdots,\xi_{k+l})=M_k(\xi_1,\cdots,
\xi_{j-1}, \xi_j+\cdots+\xi_{j+l},\xi_{j+l+1}, \cdots, \xi_{k+l}).
\endaligned
\end{equation*}
Also if $M_k$ is a multiplier of order $k$ and $u_1, u_2, \cdots,
u_k$ are functions on $\mathbb{R}^d$, we define the $k-$linear
functional
\begin{equation*}
\aligned \Lambda_k\big(M_k; u_1, u_2, \cdots u_k
\big)=\text{Re}\int_{\Gamma_k}M_k(\xi_1, \xi_2,\cdots, \xi_k)
\prod^k_{j=1}\widehat{u}_j(\xi_j)
\endaligned
\end{equation*}
and we adopt the notation $\Lambda_k\big(M_k;
u\big)=\Lambda_k\big(M_k; u, \overline{u}, \cdots, u,
\overline{u}\big)$. We observe that the quantity $\Lambda_k\big(M_k;
u\big)$ is invariant
\begin{enumerate}\label{sym}
\item[$(1)$]if one permutes the even arguments $\xi_2, \xi_4,\cdots, \xi_k$ of
$M_k$;
\item[$(2)$]if one permutes the odd arguments $\xi_1, \xi_3, \cdots, \xi_{k-1}$ of
$M_k$;
\item[$(3)$]if one makes the change of
\begin{equation*}
\aligned M_k(\xi_1, \xi_2,\cdots,\xi_{k-1}, \xi_k) \mapsto
M_k(-\xi_2, -\xi_1, \cdots, -\xi_k,-\xi_{k-1}).
\endaligned
\end{equation*}
\end{enumerate}

If $u$ is a solution of (\ref{equ1}), the following differentiation
law holds for the multiplier forms $\Lambda_k\big(M_k; u\big)$
\begin{equation}\label{differentialrule}
\aligned \partial_t \Lambda_k\big(& M_k; u\big)\\=&
\Lambda_k\big(iM_k \sum^k_{j=1} (-1)^j|\xi_j|^2; u \big) +
\Lambda_{k+2}\big(i \sum^k_{j=1}(-1)^j|\xi_{j+1,
j+2}|^{-(d-\gamma)}X^2_j(M_k) ; u \big)
\endaligned
\end{equation}
where we use the notational convention $\xi_{a,b}=\xi_a + \xi_b$,
$\xi_{a,b,c}=\xi_a + \xi_b + \xi_c$, etc.

Using the above notation, the modified energy (\ref{I1}) can be
written as follows:
\begin{equation*}
\aligned \widetilde{E}(u)=\Lambda_2\big(-\frac12\xi_1 m_1 \cdot
\xi_2 m_2; u \big) + \Lambda_4\big(\frac14|\xi_{2,3} |^{-(d-\gamma)}
m_1 m_2 m_3 m_4; u \big)
\endaligned
\end{equation*}
where we abbreviate $m(\xi_j)$ as $m_j$.

Together with the the differentiation rules (\ref{differentialrule})
and the symmetry properties of k-linear functional
$\Lambda_k\big(M_k; u\big)$, we obtain
\begin{equation*}
\aligned
\partial_t \Lambda_2\big( -\frac{1}{2} \xi_1 m_1 \cdot \xi_2 m_2; u
\big) & = \Lambda_2\big( -\frac{i}{2} \xi_1 m_1 \cdot \xi_2
m_2\sum^2_{j=1}(-1)^{j}|\xi_j|^2; u \big) \\
& \qquad + \Lambda_4\big( -\frac{i}{2}\sum^2_{j=1}(-1)^{j}
|\xi_{j+1, j+2}|^{-(d-\gamma)}X^2_j\big( \xi_1 m_1 \cdot \xi_2
m_2\big); u
\big)\\
& = \Lambda_4\big( i |\xi_{2, 3}|^{-(d-\gamma)}m^2_1 |\xi_1|^2; u
\big),
\endaligned
\end{equation*}
and
\begin{equation*}
\aligned
\partial_t \Lambda_4 \big(\frac14 |\xi_{2,3} |^{-(d-\gamma)} & m_1 m_2 m_3
m_4; u\big) \\& = \Lambda_4 \big(\frac{i}{4} |\xi_{2,3}
|^{-(d-\gamma)} m_1 m_2 m_3
m_4\sum^{4}_{j=1}(-1)^{j}|\xi_j|^2; u\big)\\
&\qquad  +  \Lambda_6 \big(\frac{i}{4}\sum^{4}_{j=1} (-1)^{j}
|\xi_{j+1, j+2}|^{-(d-\gamma)}X^2_j\big( |\xi_{2,3} |^{-(d-\gamma)}
m_1 m_2 m_3 m_4\big); u\big) \\
& =-\Lambda_4 \big(i |\xi_{2,3} |^{-(d-\gamma)}|\xi_1|^2 m_1 m_2 m_3
m_4; u\big)\\
&\qquad  -  \Lambda_6 \big(i|\xi_{2, 3}|^{-(d-\gamma)} |\xi_{4,
5}|^{-(d-\gamma)}  m_{1,2,3} m_4 m_5 m_6; u\big)\\
& =-\Lambda_4 \big(i |\xi_{2,3} |^{-(d-\gamma)}|\xi_1|^2 m_1 m_2 m_3
m_4; u\big)\\
&\qquad  + \Lambda_6 \big(i|\xi_{2, 3}|^{-(d-\gamma)} |\xi_{4,
5}|^{-(d-\gamma)}  m_{1,2,3}\big( m_{1,2,3}- m_4 m_5 m_6\big);
u\big).
\endaligned
\end{equation*}

The fundamental theorem of calculus  together with these estimates
implies the following proposition, which will be used to prove that
$\widetilde{E}$ is almost conserved.

\begin{proposition}\label{increformular}
Let $u$ be an $H^1$ solution to (\ref{equ1}). Then for any $T\in
\mathbb{R}$ and $\delta>0$, we have
\begin{equation*}
\aligned
\widetilde{E}(u)(T+\delta)-\widetilde{E}(u)(T)=\int^{T+\delta}_{T}\Lambda_4\big(M_4;
u \big)\ dt+ \int^{T+\delta}_{T}\Lambda_6\big(M_6; u \big)\ dt
\endaligned
\end{equation*}
with \begin{equation*} \aligned M_4&=i \big|\xi_{2,
3}\big|^{-(d-\gamma)} \big|\xi_1\big|^2m_1\big(m_1- m_2 m_3
m_4\big);\\
M_6&=i\big|\xi_{2, 3}\big|^{-(d-\gamma)} \big|\xi_{4,
5}\big|^{-(d-\gamma)} m_{1,2,3}\big( m_{1,2,3}- m_4 m_5 m_6\big).
\endaligned
\end{equation*}
Furthermore if $|\xi_j| \ll N$ for all $j$, then the multipliers
$M_4$ and $M_6$ vanish on $\Gamma_4$ and $\Gamma_6$, respectively.
\end{proposition}

\subsection{Modified local well-posedness}
In this subsection, we shall prove a local well-posedness result for
the modified solution $Iu$ and some a priori estimates for it.

Let $J=[t_0, t_1]$ be an interval of time. We denote by $Z_I(J)$ the
following space:
\begin{equation*}
\aligned Z_I(J)=S_I(J)\cap X^{1,\frac12+}_I(J)
\endaligned
\end{equation*}
where
\begin{equation*}
\aligned S_I(J)&=\Big\{u;  \sup_{(q, r) \ \text{admissible}} \big\| \langle \nabla \rangle Iu\big\|_{L^q_tL^r_x(J\times \mathbb{R}^d)}< \infty \Big\}, \\
&X^{1,b}_I(J) = \Big\{u; \quad
\big\|Iu\big\|_{X^{1,\frac12+}(J\times \mathbb{R}^d)}< \infty
\Big\}.
\endaligned
\end{equation*}

\begin{proposition}\label{MLWP}
Let $2<\gamma<3 \leq d$, $s>\frac{\gamma}{2}-1$, and consider the
IVP
\begin{equation}\label{mequ}
\aligned
iIu_t+\Delta Iu &=I\big(|x|^{-\gamma}*|u|^2 u \big),\ x\in \mathbb{R}^d, \ t\in \mathbb{R}\\
Iu(t_0, x)&=Iu_0(x)\in H^1(\mathbb{R}^d).
\endaligned
\end{equation}
Then for any $u_0\in H^s$, there exists a time interval $J=[t_0,
t_0+\delta]$, $\delta=\delta(\big\|Iu_0\big\|_{H^1})$ and there
exists a unique $u\in Z_I(J)$ solution to (\ref{mequ}). Moreover
there is continuity with respect to the initial data.
\end{proposition}

{\bf Proof: } The proof of this proposition proceeds by the usual
fixed point method on the space $Z_I(J)$. Since the estimates are
very similar to the ones we provide in the proof of Proposition
\ref{universalbound} below, in particular (\ref{ine5}) and
(\ref{ine6}) , we omit the details.

\begin{proposition}\label{universalbound}
Let $2<\gamma<3\leq d$ and $s>\frac{\gamma}{2}-1$. If $u$ is a
solution to the IVP (\ref{mequ}) on the interval $J=[t_0, t_1]$,
which satisfies the following a priori bound
\begin{equation*}
\aligned \big\|I
u\big\|^{4}_{L^{4}_t\dot{H}^{-\frac{d-3}{4},4}_x\big(J\times
\mathbb{R}^d \big)} < \mu,
\endaligned
\end{equation*}
where $\mu$ is a small universal constant, then
\begin{equation*}
\aligned \big\|u\big\|_{Z_I(J)}\lesssim\big\|Iu_0\big\|_{H^1}.
\endaligned
\end{equation*}
\end{proposition}

{\bf Proof: } We start by obtaining a control of the Strichartz
norms. Applying $\langle \nabla \rangle$ to (\ref{mequ}) and using
the Strichartz estimate in Proposition \ref{stri}. For any pair of
admissible exponents $(q,r)$, we obtain
\begin{equation}\label{ine1}
\aligned \big\| \langle \nabla \rangle I u\big\|_{L^q_tL^r_x}
\lesssim \big\|Iu_0 \big\|_{H^1} + \big\| \langle \nabla \rangle I
\big( (|x|^{-\gamma}*|u|^2)u\big)
\big\|_{L^{\frac32}_tL^{\frac{6d}{3d+4}}_x}.
\endaligned
\end{equation}
Now we notice that the multiplier $ \langle \nabla \rangle I $ has
symbol which is increasing as a function of $|\xi|$ for any $s\geq
\frac{\gamma}{2}-1$. Using this fact one can modify the proof of the
Leibnitz rule for fractional derivatives and prove its validity for
$ \langle \nabla \rangle I$. See also Principle A.5 in the appendix
of \cite{Taobook1}. This remark combined with (\ref{ine1}) implies
that
\begin{equation}\label{ine4}
\aligned \big\| \langle \nabla \rangle I u\big\|_{L^q_tL^r_x} &
\lesssim \big\|Iu_0 \big\|_{H^1} + \big\| \langle \nabla \rangle I
\big( (|x|^{-\gamma}*|u|^2)u\big)
\big\|_{L^{\frac32}_tL^{\frac{6d}{3d+4}}_x} \\
& \lesssim \big\|Iu_0 \big\|_{H^1} + \big\||x|^{-\gamma}* \langle
\nabla \rangle I (|u|^2) \big\|_{L^2_tL^{\frac{2d}{\gamma}}_x}
\big\| u \big\|_{L^6_tL^{\frac{6d}{3d+4-3\gamma}}_x}  \\
& \qquad \qquad \quad + \big\||x|^{-\gamma}* |u|^2
\big\|_{L^3_tL^{\frac{3d}{4}}_x} \big\| \langle \nabla \rangle I u
\big\|_{L^3_tL^{\frac{6d}{3d-4}}_x} \\
& \lesssim \big\|Iu_0 \big\|_{H^1} +  \big\| \langle \nabla \rangle
I u \big\|_{L^3_tL^{\frac{6d}{3d-4}}_x} \big\| u
\big\|^2_{L^6_tL^{\frac{6d}{3d+4-3\gamma}}_x}
\endaligned
\end{equation}
where we used H\"{o}lder's inequality and Hardy-Littlewood-Sobolev's
inequaltiy.

In order to obtain an upper bound on $ \big\| u
\big\|_{L^6_tL^{\frac{6d}{3d+4-3\gamma}}_x}$, we perform a
Littlewood-Paley decomposition along the following lines. We note
that a similar approach was used in \cite{CKSTT04}. We write
\begin{equation}\label{decomp}
\aligned u=u_{N_0} + \sum^{\infty}_{j=1} u_{N_j},
\endaligned
\end{equation}
where $u_{N_0}$ has spatial frequency support for $\langle \xi
\rangle \leq N$, while $u_{N_j}$ is such that its spatial Fourier
support transform is supported for $\langle \xi \rangle \thickapprox
N_j =2^{h_j}$ with $h_j \gtrsim \log N$ and $j=1,2, \cdots$. By the
triangle inequality and H\"{o}lder's inequality, we have
\begin{equation}\label{ine2}
\aligned \big\| u \big\|_{L^6_tL^{\frac{6d}{3d+4-3\gamma}}_x}
\lesssim & \big\| u_{N_0}
\big\|_{L^6_tL^{\frac{6d}{3d+4-3\gamma}}_x}+ \sum^{\infty}_{j=1}
\big\| u_{N_j}
\big\|_{L^6_tL^{\frac{6d}{3d+4-3\gamma}}_x} \\
\lesssim &  \big\| u_{N_0}
\big\|_{L^6_tL^{\frac{6d}{3d+4-3\gamma}}_x} + \sum^{\infty}_{j=1}
\big\| u_{N_j}
\big\|^{2-\frac{\gamma}{2}}_{L^6_tL^{\frac{6d}{3d-2}}_x}\big\|
u_{N_j} \big\|^{\frac{\gamma}{2}-1}_{L^6_tL^{\frac{6d}{3d-8}}_x}.
\endaligned
\end{equation}
On the other hand, by using the definition of the operator $I$, the
defintion of the $u_{N_j}$'s and the Marcinkiewicz multiplier
theorem, we observe that for some $0<\theta_i<1, i=1,\cdots, 4,
\sum^4_{j=1}\theta_i=1$

\begin{equation*}
\aligned  \big\| u_{N_0} \big\|_{L^6_tL^{\frac{6d}{3d+4-3\gamma}}_x}
& \lesssim  \big\|
u_{N_0}\big\|^{\theta_1}_{L^{4}_t\dot{H}^{-\frac{d-3}{4},4}_x} \big\|u_{N_0}\big\|^{\theta_2}_{L^6_tL^{\frac{6d}{3d-8}}_x} \big\|u_{N_0}\big\|^{\theta_3}_{L^6_tL^{\frac{6d}{3d-2}}_x}\big\|u_{N_0}\big\|^{\theta_4}_{L^{\infty}_tL^2_x}\\
&  \lesssim  \big\|I
u_{N_0}\big\|^{\theta_1}_{L^{4}_t\dot{H}^{-\frac{d-3}{4},4}_x}
\big\| u_{N_0}\big\|^{1-\theta_1}_{Z_I(J)}.\\
\big\| \langle \nabla \rangle I u_{N_j}
\big\|_{L^6_tL^{\frac{6d}{3d-2}}_x} & \thickapprox N_j
\big(\frac{N}{N_j} \big)^{1-s} \big\| u_{N_j}
\big\|_{L^6_tL^{\frac{6d}{3d-2}}_x},
\quad j=1,2,\cdots.\\
 \big\| I u_{N_j}
\big\|_{L^6_tL^{\frac{6d}{3d-8}}_x} & \thickapprox\big(\frac{N}{N_j}
\big)^{1-s}\big\| u_{N_j} \big\|_{L^6_tL^{\frac{6d}{3d-8}}_x}, j=1,
2,\cdots.
\endaligned
\end{equation*}

Now we use these estimates to obtain the following upper bound on
(\ref{ine2})

\begin{equation}\label{ine3}
\aligned \big\| u \big\|_{L^6_tL^{\frac{6d}{3d+4-3\gamma}}_x}
\lesssim &\  \big\|I
u_{N_0}\big\|^{\theta_1}_{L^{4}_t\dot{H}^{-\frac{d-3}{4},4}_x}
\big\|u_{N_0}\big\|^{1-\theta_1}_{Z_I(J)} \\
& \ \ + \sum^{\infty}_{j=1}\Big( \frac{1}{N_j}
\big(\frac{N_j}{N}\big)^{1-s}\big\|u_{N_j}\big\|_{Z_I(J)}
\Big)^{2-\frac{\gamma}{2}}\Big( \big(
\frac{N_j}{N}\big)^{1-2}\big\|u_{N_j}\big\|_{Z_I(J)}\Big)^{\frac{\gamma}{2}-1}\\
\lesssim &\ \mu^{\frac{\theta_1}{4}}
\big\|u\big\|^{1-\theta_1}_{Z_I(J)}+N^{-(2-\frac{\gamma}{2})}\big\|u
\big\|_{Z_I(J)},
\endaligned
\end{equation}
which together with (\ref{ine4}) implies that
\begin{equation}\label{ine5}
\aligned \big\| \langle \nabla \rangle I u\big\|_{L^q_tL^r_x} &
\lesssim \big\|Iu_0 \big\|_{H^1} + \mu^{\frac{\theta_1}{2}}\big\|u
\big\|^{3-2\theta_1}_{Z_I(J)}+N^{-(4-\gamma)}\big\|u
\big\|^3_{Z_I(J)}.
\endaligned
\end{equation}

Now we shall obtain a control of the $X^{s,b}$ norm. We use
Duhamel's formula and the theory of $X^{s,b}$ spaces \cite{Grun02},
\cite{Taobook1} to obtain
\begin{equation}\label{ine6}
\aligned \big\|Iu\big\|_{X^{1,\frac12+}} & \lesssim \big\|Iu_0
\big\|_{H^1} +\big\| \langle \nabla \rangle I \big(
(|x|^{-\gamma}*|u|^2)u\big)
\big\|_{X^{0,-\frac{1}{2}+}}\\
& \lesssim \big\|Iu_0 \big\|_{H^1}+ \big\| \langle \nabla \rangle I
\big( (|x|^{-\gamma}*|u|^2)u\big)
\big\|_{L^{\frac32+}_tL^{\frac{6d}{3d+4}+}_x}\\
& \lesssim \big\|Iu_0 \big\|_{H^1} +  \big\| \langle \nabla \rangle
I u \big\|_{L^3_tL^{\frac{6d}{3d-4}}_x} \big\| u
\big\|_{L^6_tL^{\frac{6d}{3d+4-3\gamma}}_x}\big\| u
\big\|_{L^{6+}_tL^{\frac{6d}{3d+4-3\gamma}+}_x}.
\endaligned
\end{equation}
An upper bound on $ \big\| u
\big\|_{L^6_tL^{\frac{6d}{3d+4-3\gamma}}_x}$ is given by
(\ref{ine3}). In order to obtain an upper bound on $\big\| u
\big\|_{L^{6+}_tL^{\frac{6d}{3d+4-3\gamma}+}_x}$, we proceed as
follows. First we perform a dyadic decomposition and write $u$ as
(\ref{decomp}). The triangle inequality applied on (\ref{decomp})
gives for any $0<\delta<\frac{\gamma}{2}-1$
\begin{equation}\label{ine7}
\aligned \big\| u & \big\|_{L^{6+}_t
L^{\frac{6d}{3d+4-3\gamma}+}_x}
\\ &\lesssim  \big\| u_{N_0}
\big\|_{L^{6+}_tL^{\frac{6d}{3d+4-3\gamma}+}_x} +
\sum^{\infty}_{j=1}\big\| u_{N_j}
\big\|_{L^{6+}_tL^{\frac{6d}{3d+4-3\gamma}+}_x} \\
&= \big\| I u_{N_0} \big\|_{L^{6+}_tL^{\frac{6d}{3d+4-3\gamma}+}_x}
+ \sum^{\infty}_{j=1} N^{\delta-s}_{j}N^{s-1}\big\|\langle \nabla
\rangle^{1-\delta}I u_{N_j}
\big\|_{L^{6+}_tL^{\frac{6d}{3d+4-3\gamma}+}_x}\\
& \lesssim  \big\| I u
\big\|_{L^{6+}_tL^{\frac{6d}{3d+4-3\gamma}+}_x} + \big\|\langle
\nabla \rangle^{1-\delta}I u
\big\|_{L^{6+}_tL^{\frac{6d}{3d+4-3\gamma}+}_x} \lesssim  \big\| I u
\big\|_{X^{1,\frac12+}},
\endaligned
\end{equation}
where we use Proposition \ref{linear}. By applying the inequalities
(\ref{ine3}) and (\ref{ine7}) to bound the right hand side of
(\ref{ine6}), we obtain
\begin{equation}\label{ine8}
\aligned \big\|Iu\big\|_{X^{1,\frac12+}} & \lesssim \big\|Iu_0
\big\|_{H^1} +\mu^{\frac{\theta_1}{4}}\big\|u
\big\|^{3-\theta_1}_{Z_I(J)}+N^{-(2-\frac{\gamma}{2})}\big\|u
\big\|^3_{Z_I(J)}.
\endaligned
\end{equation}
The desired bound follows from (\ref{ine5}) and (\ref{ine8}) by
choosing $N$ sufficiently large.

\subsection{An upper bound on the increment of $\widetilde{E}(u)$}
{\bf Decomposition remark. } Our approach to prove a decay for the
increment of the modified energy is based on obtaining certain
multilinear estimates in appropriate functional spaces which are
$L^2$-based. Hence, whenever we perform a Littlewood-Paley
decomposition of a function we shall assume that the Fourier
transforms of the Littlewood-Paley pieces are positive. Moreover, we
will ignore the presence of conjugates. At the end we will always
keep a decay factor $C(N_1, N_2, \cdots)$ in order to perform the
summations.

Now we proceed to prove the almost conservation law of the modified
energy. In Proposition \ref{increformular}, we prove that an
increment of the modified energy can be expressed as
\begin{equation*}
\aligned
\widetilde{E}(u)(T+\delta)-\widetilde{E}(u)(T)=\int^{T+\delta}_{T}\Lambda_4\big(M_4;
u \big)\ dt+ \int^{T+\delta}_{T}\Lambda_6\big(M_6; u \big)\ dt
\endaligned
\end{equation*}
with \begin{equation*} \aligned M_4&=i |\xi_{2, 3}|^{-(d-\gamma)}
|\xi_1|^2m_1\big( m_1- m_2 m_3
m_4\big);\\
M_6&=i|\xi_{2, 3}|^{-(d-\gamma)} |\xi_{4, 5}|^{-(d-\gamma)}
m_{1,2,3}\big( m_{1,2,3}- m_4 m_5 m_6\big).
\endaligned
\end{equation*}
Hence in order to control the increment of the modified energy, we
shall find an upper bound on the $\Lambda_4\big(M_4; u \big)$ and
$\Lambda_6\big(M_6; u \big)$ forms, which we do in the following
propositions. First we give the quadrilinear  estimate.

\begin{proposition}\label{increment1}
For any Schwartz function $u$, and any $\delta \thickapprox 1$ just
as in Proposition \ref{MLWP}, we have that
\begin{equation}\label{bound1}
\aligned \Big| \int^{T+\delta}_{T}\Lambda_4\big(M_4; u \big)\ dt
\Big| \lesssim N^{-1+} \big\|Iu \big\|^4_{X^{1, \frac12+}},
\endaligned
\end{equation}
for $s>\frac{\gamma}{2}-1$.
\end{proposition}

{\bf Proof: } By Plancherel theorem, we aim to prove that
\begin{equation}\label{goal1}
\aligned \Big| \int^{T+\delta}_{T} \int_{\Gamma_4}|\xi_{2,
3}|^{-(d-\gamma)}|\xi_1|^2 m_1\big( m_1- &m_2 m_3  m_4\big)
\widehat{u}_1(t,\xi_1) \widehat{\overline{u}}_2(t,\xi_2)
\widehat{u}_3(t,\xi_3)\widehat{\overline{u}}_4(t,\xi_4) \Big| \\
\lesssim & \ N^{-1+} C(N_1, N_2, N_3, N_4) \prod^4_{j=1}\big\|Iu_j
\big\|_{X^{1, \frac12+}},
\endaligned
\end{equation}
where $C(N_1, N_2, N_3, N_4)$ is a decay just as the remark above,
and it allows us to sum over all dyadic shells. The analysis which
follows will not rely on the complex conjugate structure in
$\Lambda_4\big(M_4; u \big)$. Thus, by symmetry, we may assume that
$N_2 \geq N_3 \geq N_4$.

{\bf Case 1: $N \gg N_2$.} According to the definition of $m(\xi)$,
the multiplier
\begin{equation*}
\aligned |\xi_{2, 3}|^{-(d-\gamma)} m_1\big( m_1- m_2 m_3 m_4\big)
\endaligned
\end{equation*}
is identically 0, the bound (\ref{bound1}) holds trivially.

{\bf Case 2: $N_2 \gtrsim N \gg N_3 \geq N_4$.} Since $\displaystyle
\sum^4_{j=1}\xi_j=0$, we have $N_1 \thickapprox N_2$. We aim for
(\ref{goal1}) with a decay factor
\begin{equation*}
\aligned C(N_1, N_2, N_3, N_4)=N^{0-}_2.
\endaligned
\end{equation*}

By the mean value theorem, we have the following pointwise bound
\begin{equation*}
\aligned \big| m_1\big( m_1- m_2 m_3 m_4\big) \big|& = \big|
m_1\big( m_{2,3,4}- m_2 m_3 m_4\big) \big|\\
& \lesssim m_1 |\nabla m(\xi)\cdot (\xi_3 + \xi_4)| \quad
\text{where} \quad
|\xi|\backsim |\xi_2|\\
& \lesssim m_1 m_2 \frac{N_3}{N_2}.
\endaligned
\end{equation*}
Hence by H\"{o}lder's inequality and Hardy-Littlewood-Sobolev's
inequality and Proposition \ref{linear}, we obtain
\begin{equation*}
\aligned \text{LHS of}\  (\ref{goal1}) & \lesssim N^2_1 m_1 m_2
\frac{N_3}{N_2} \Big|\int^{T+\delta}_{T} \int_{\Gamma_4}  |\xi_{2,
3}|^{-(d-\gamma)} \widehat{u}_1(t,\xi_1)
\widehat{\overline{u}}_2(t,\xi_2)
\widehat{u}_3(t,\xi_3)\widehat{\overline{u}}_4(t,\xi_4)\Big| \\
& \lesssim N^2_1 m_1 m_2 \frac{N_3}{N_2} \big\|u_1
\big\|_{L^{3}_tL^{\frac{6d}{3d-4}}_x} \big\|u_2
\big\|_{L^{3}_tL^{\frac{6d}{3d-4}}_x} \big\|u_3
\big\|_{L^{6}_tL^{\frac{6d}{3d-2}}_x} \big\|u_4
\big\|_{L^{6}_tL^{\frac{6d}{3d+10-6\gamma}}_x} \\
& \lesssim N^2_1 m_1 m_2 \frac{N_3}{N_2} N^{\gamma-2}_4
\prod^4_{j=1}\big\|u_j \big\|_{X^{0, \frac12+}}.
\endaligned
\end{equation*}

It suffices to show that
\begin{equation*}
\aligned N^2_1 m_1 m_2 \frac{N_3}{N_2} N^{\gamma-2}_4 \lesssim
N^{-1+} N^{0-}_2 \ m_1 N_1 m_2 N_2 \langle N_3 \rangle \langle N_4
\rangle.
\endaligned
\end{equation*}
We reduce to show that
\begin{equation*}
\aligned N^{1-}N^{0+}_2 \lesssim N_2 \ \langle N_3\rangle N^{-1}_3\
\langle N_4\rangle N^{2-\gamma}_4.
\endaligned
\end{equation*}
This is true since
\begin{equation*}
\aligned N_2  \gtrsim & N^{1-}N^{0+}_2 ;\\
\langle N_3\rangle N^{-1}_3  \gtrsim 1; &\quad  \langle N_4\rangle
N^{2-\gamma}_4  \gtrsim 1.
\endaligned
\end{equation*}

{\bf Case 3: $N_2 \geq N_3 \gtrsim N $.} In this case, we use the
trivial pointwise bound
\begin{equation*}
\aligned \big| m_1\big( m_1- m_2 m_3 m_4\big) \big|& \lesssim m^2_1.
\endaligned
\end{equation*}
The frequency interactions fall into two subcategories, depending on
which frequency is comparable to $N_2$.

{\bf Case 3a: $N_1 \thickapprox N_2 \geq N_3 \gtrsim N $.} In this
case, we prove the decay factor
\begin{equation*}
\aligned
 C(N_1, N_2, N_3, N_4)=N^{0-}_3.
\endaligned
\end{equation*}
in (\ref{goal1}). This allows us to directly sum in $N_3$ and $N_4$,
and sum in $N_1$ and $N_2$ after applying Cauchy-Schwarz to those
factors.

By H\"{o}lder's inequality and Hardy-Littlewood-Sobolev's inequality
and Proposition \ref{linear}, we obtain
\begin{equation*}
\aligned \text{LHS of}\  (\ref{goal1}) & \lesssim N^2_1 m^2_1
\Big|\int^{T+\delta}_{T} \int_{\Gamma_4}  |\xi_{2, 3}|^{-(d-\gamma)}
\widehat{u}_1(t,\xi_1) \widehat{\overline{u}}_2(t,\xi_2)
\widehat{u}_3(t,\xi_3)\widehat{\overline{u}}_4(t,\xi_4)\Big| \\
& \lesssim  N^2_1 m^2_1 \big\|u_1
\big\|_{L^{3}_tL^{\frac{6d}{3d-4}}_x} \big\|u_2
\big\|_{L^{3}_tL^{\frac{6d}{3d-4}}_x} \big\|u_3
\big\|_{L^{6}_tL^{\frac{6d}{3d-2}}_x} \big\|u_4
\big\|_{L^{6}_tL^{\frac{6d}{3d+10-6\gamma}}_x} \\
& \lesssim  N^2_1 m^2_1 N^{\gamma-2}_4 \prod^4_{j=1}\big\|u_j
\big\|_{X^{0, \frac12+}}.
\endaligned
\end{equation*}

It suffices to show that
\begin{equation*}
\aligned N^2_1 m^2_1 N^{\gamma-2}_4 \lesssim N^{-1+} N^{0-}_3 \ m_1
N_1 m_2 N_2  m_3 N_3 \langle N_4 \rangle.
\endaligned
\end{equation*}
We reduce to show that
\begin{equation*}
\aligned N^{1-}N^{0+}_3 \lesssim m_3  N_3\ m_4\langle N_4\rangle
N^{2-\gamma}_4.
\endaligned
\end{equation*}
This is true since for $s\geq \gamma-2$, we have
\begin{equation*}
\aligned m_3  N_3\ m_4\langle N_4\rangle N^{2-\gamma}_4 & \geq m_3
N_3\ m_4\langle N_4\rangle^{3-\gamma} \\
& \gtrsim m_3 N_3 \gtrsim N^{1-}N^{0+}_3,
\endaligned
\end{equation*}
where we used the fact that $m(\xi)\langle\xi\rangle^p$ is monotone
non-decreasing if $s+p\geq 1$. While for $\frac{\gamma}{2}-1 < s <
\gamma-2$, we have
\begin{equation*}
\aligned m_3  N_3\ m_4\langle N_4\rangle N^{2-\gamma}_4 & \gtrsim
m_3  N_3\ m_3 N^{3-\gamma}_3\\
& \gtrsim N^{4-\gamma-}N^{0+}_3 \gtrsim N^{1-}N^{0+}_3,
\endaligned
\end{equation*}
where we used the fact that $m(\xi)\langle\xi\rangle^p$ is monotone
non-increasing if $s+p < 1$.

{\bf Case 3b: $ N_2 \thickapprox N_3 \gtrsim N , N_2\gtrsim N_1$.}
In this case, we prove the decay factor
\begin{equation*}
\aligned
 C(N_1, N_2, N_3, N_4)=N^{0-}_2.
\endaligned
\end{equation*}
in (\ref{goal1}). This will allow us to directly sum in all the
$N_i$.

By H\"{o}lder's inequality and Hardy-Littlewood-Sobolev's inequality
and Proposition \ref{linear} once again, we obtain
\begin{equation*}
\aligned \text{LHS of}\  (\ref{goal1}) & \lesssim N^2_1 m^2_1
\Big|\int^{T+\delta}_{T} \int_{\Gamma_4}  |\xi_{2, 3}|^{-(d-\gamma)}
\widehat{u}_1(t,\xi_1) \widehat{\overline{u}}_2(t,\xi_2)
\widehat{u}_3(t,\xi_3)\widehat{\overline{u}}_4(t,\xi_4)\Big| \\
& \lesssim  N^2_1 m^2_1 \big\|u_1
\big\|_{L^{3}_tL^{\frac{6d}{3d-4}}_x} \big\|u_2
\big\|_{L^{3}_tL^{\frac{6d}{3d-4}}_x} \big\|u_3
\big\|_{L^{6}_tL^{\frac{6d}{3d-2}}_x} \big\|u_4
\big\|_{L^{6}_tL^{\frac{6d}{3d+10-6\gamma}}_x} \\
& \lesssim  N^2_1 m^2_1 N^{\gamma-2}_4 \prod^4_{j=1}\big\|u_j
\big\|_{X^{0, \frac12+}}.
\endaligned
\end{equation*}

It suffices to show that
\begin{equation*}
\aligned N^2_1 m^2_1 N^{\gamma-2}_4 \lesssim N^{-1+} N^{0-}_2 \ m_1
N_1 m_2 N_2  m_3 N_3 \langle N_4 \rangle.
\endaligned
\end{equation*}
Note that
\begin{equation*}
\aligned N^2_1 m^2_1 \lesssim m_1 N_1 m_2 N_2.
\endaligned
\end{equation*}
We reduce to show that
\begin{equation*}
\aligned N^{1-}N^{0+}_2 \lesssim m_3  N_3\ m_4\langle N_4\rangle
N^{2-\gamma}_4.
\endaligned
\end{equation*}
This is true since for $s\geq \gamma-2$, we have
\begin{equation*}
\aligned m_3  N_3\ m_4\langle N_4\rangle N^{2-\gamma}_4 & \geq m_3
N_3\ m_4\langle N_4\rangle^{3-\gamma} \\
& \gtrsim m_3 N_3 \approx m_2 N_2 \gtrsim N^{1-}N^{0+}_2,
\endaligned
\end{equation*}
where we used the fact that $m(\xi)\langle\xi\rangle^p$ is monotone
non-decreasing if $s+p\geq 1$. While for $\frac{\gamma}{2}-1 < s <
\gamma-2$, we have
\begin{equation*}
\aligned m_3  N_3\ m_4\langle N_4\rangle N^{2-\gamma}_4 & \gtrsim
m_3  N_3\ m_3 N^{3-\gamma}_3\approx m^2_2 N^{4-\gamma}_2\\
& \gtrsim N^{4-\gamma-}N^{0+}_3 \gtrsim N^{1-}N^{0+}_3,
\endaligned
\end{equation*}
where we used the fact that $m(\xi)\langle\xi\rangle^p$ is monotone
non-increasing if $s+p < 1$. This completes the proof.

In order to make use of quadrilinar estimate (Proposition
\ref{increment1}) to obtain sextilinear estimate, we first give a
lemma
\begin{lemma}\label{trilinear}
Assume $u$, $\delta$ are as in Proposition \ref{MLWP}, and
$P_{N_{1,2,3}}$ the Littlewood-Paley projection onto the $N_{1,2,3}$
frequency shell. Then
\begin{equation*}
\aligned \big\|P_{N_{1,2,3}}\big( I \big(u
|\nabla|^{-(d-\gamma)}|u|^2\big)\big)\big\|_{L^3_tL^{\frac{6d}{3d-4}}_x}
\lesssim N_{1,2,3}\big\|Iu\big\|^3_{X^{1, \frac12+}}.
\endaligned
\end{equation*}
\end{lemma}

{\bf Proof: } We write $u=u_L+u_H$ where
\begin{equation*}
\aligned \text{supp} \ \widehat{u}_l(t,\xi) & \subseteq \big\{ |\xi|<2 \big\}, \\
\text{supp} \ \widehat{u}_H(t,\xi) & \subseteq \big\{ |\xi|>1
\big\}.
\endaligned
\end{equation*}
Hence,
\begin{equation*}
\aligned  &\big\| P_{N_{1,2,3}} \big( I \big(u
|\nabla|^{-(d-\gamma)}|u|^2\big)\big)\big\|_{L^3_tL^{\frac{6d}{3d-4}}_x}\\
 \lesssim & \quad \ \big\|P_{N_{1,2,3}}\big( I \big(u_L
|\nabla|^{-(d-\gamma)}|u_L|^2\big)\big)\big\|_{L^3_tL^{\frac{6d}{3d-4}}_x}
 + \big\|P_{N_{1,2,3}}\big( I \big(u_H
|\nabla|^{-(d-\gamma)}|u_H|^2\big)\big)\big\|_{L^3_tL^{\frac{6d}{3d-4}}_x}
\\
& + \big\|P_{N_{1,2,3}}\big( I \big(u_L
|\nabla|^{-(d-\gamma)}|u_H|^2\big)\big)\big\|_{L^3_tL^{\frac{6d}{3d-4}}_x}
+ \big\|P_{N_{1,2,3}}\big( I \big(u_H
|\nabla|^{-(d-\gamma)}|u_L|^2\big)\big)\big\|_{L^3_tL^{\frac{6d}{3d-4}}_x}\\
& + \big\|P_{N_{1,2,3}}\big( I \big(u_L
|\nabla|^{-(d-\gamma)}\overline{u}_Lu_H\big)\big)\big\|_{L^3_tL^{\frac{6d}{3d-4}}_x}
+ \big\|P_{N_{1,2,3}}\big( I \big(u_L
|\nabla|^{-(d-\gamma)}u_L\overline{u}_H\big)\big)\big\|_{L^3_tL^{\frac{6d}{3d-4}}_x}\\
&+ \big\|P_{N_{1,2,3}}\big( I \big(u_H
|\nabla|^{-(d-\gamma)}\overline{u}_Lu_H\big)\big)\big\|_{L^3_tL^{\frac{6d}{3d-4}}_x}
+ \big\|P_{N_{1,2,3}}\big( I \big(u_H
|\nabla|^{-(d-\gamma)}u_L\overline{u}_H\big)\big)\big\|_{L^3_tL^{\frac{6d}{3d-4}}_x}
\endaligned
\end{equation*}

Consider the first term.  By H\"{o}lder's inequality,
Hardy-Littlewood-Sobolev's inequality and Proposition \ref{linear},
we have
\begin{equation*}
\aligned \big\|P_{N_{1,2,3}}\big( I \big(u_L
|\nabla|^{-(d-\gamma)}|u_L|^2\big)\big)\big\|_{L^3_tL^{\frac{6d}{3d-4}}_x}
& \lesssim \big\|u_L
|\nabla|^{-(d-\gamma)}|u_L|^2\big\|_{L^3_tL^{\frac{6d}{3d-4}}_x} \\
& \lesssim \big\|u_L\big\|^3_{L^9_tL^{\frac{18d}{9d-6\gamma-4}}_x}=
\big\|I u_L\big\|^3_{L^9_tL^{\frac{18d}{9d-6\gamma-4}}_x} \\
& \lesssim \big\|I u_L\big\|^3_{X^{1,\frac12+}} \leq
N_{1,2,3}\big\|Iu\big\|^3_{X^{1, \frac12+}}
\endaligned
\end{equation*}
since \begin{equation*} \aligned N_{1,2,3}\geq 1, \quad
\text{and}\quad 0 \leq d \times
\big(\frac12-\frac{9d-6\gamma-4}{18d}\big)-\frac{2}{9}
=\frac{\gamma}{3} \leq 1.
\endaligned
\end{equation*}

We estimate the second term. By Sobolev's inequality and using the
Leibniz rule for the operator $|\nabla|^{2-\frac{\gamma}{2}}  I$ and
Proposition \ref{linear}, we have
\begin{equation*}
\aligned \big\|\frac{1}{N_{1,2,3}}P_{N_{1,2,3}}\big( I & \big(u_H
|\nabla|^{-(d-\gamma)}|u_H|^2\big)\big)\big\|_{L^3_tL^{\frac{6d}{3d-4}}_x}
\\ & \lesssim  \big\||\nabla|^{-1}P_{N_{1,2,3}}\big( I \big(u_H
|\nabla|^{-(d-\gamma)}|u_H|^2\big)\big)\big\|_{L^3_tL^{\frac{6d}{3d-4}}_x}
\\
& \lesssim  \big\||\nabla|^{2-\frac{\gamma}{2}}  I \big(u_H
|\nabla|^{-(d-\gamma)}|u_H|^2
\big)\big\|_{L^3_tL^{\frac{6d}{3d-3\gamma+14}}_x} \\
& \lesssim \big\||\nabla|^{2-\frac{\gamma}{2}}  I
u_H\big\|_{L^9_tL^{\frac{18d}{9d-9\gamma+14}}_x}\big\|
u_H\big\|^2_{L^9_tL^{\frac{18d}{9d-9\gamma+14}}_x}
\\
& \lesssim \big\||\nabla|^{2-\frac{\gamma}{2}}  I
u_H\big\|^3_{L^9_tL^{\frac{18d}{9d-9\gamma+14}}_x}\lesssim
\big\||\nabla|I u_H\big\|^3_{L^9_tL^{\frac{18d}{9d-4}}_x} \leq
\big\|Iu\big\|^3_{X^{1, \frac12+}}.
\endaligned
\end{equation*}

 As for the third term. By Sobolev's inequality and using the
Leibniz rule for the operator $|\nabla|^{2-\frac{\gamma}{2}}  I$ and
Proposition \ref{linear} again, we have
\begin{equation*}
\aligned \big\|\frac{1}{N_{1,2,3}}P_{N_{1,2,3}} \big( I  & \big(u_L
|\nabla|^{-(d-\gamma)}|u_H|^2\big)\big)\big\|_{L^3_tL^{\frac{6d}{3d-4}}_x}\\
& \lesssim  \big\||\nabla|^{-1}P_{N_{1,2,3}}\big( I \big(u_L
|\nabla|^{-(d-\gamma)}|u_H|^2\big)\big)\big\|_{L^3_tL^{\frac{6d}{3d-4}}_x}
\\
& \lesssim  \big\||\nabla|^{2-\frac{\gamma}{2}}  I \big(u_L
|\nabla|^{-(d-\gamma)}|u_H|^2
\big)\big\|_{L^3_tL^{\frac{6d}{3d-3\gamma+14}}_x} \\
& \lesssim \big\||\nabla|^{2-\frac{\gamma}{2}}  I
u_H\big\|_{L^9_tL^{\frac{18d}{9d-9\gamma+14}}_x}\big\|u_H
\big\|_{L^9_tL^{\frac{18d}{9d-4}}_x}\big\|u_L
\big\|_{L^9_tL^{\frac{18d}{9d-18\gamma+32}}_x} \\
& \qquad + \big\|u_H\big\|^2_{L^9_tL^{\frac{18d}{9d-9\gamma+14}}_x}
\big\||\nabla|^{2-\frac{\gamma}{2}}
Iu_L\big\|_{L^9_tL^{\frac{18d}{9d-9\gamma+14}}_x}\\
& \lesssim \big\||\nabla|I
u_H\big\|_{L^9_tL^{\frac{18d}{9d-4}}_x}\big\||\nabla|Iu_H
\big\|_{L^9_tL^{\frac{18d}{9d-4}}_x}\big\||\nabla|^{\gamma-2}Iu_L
\big\|_{L^9_tL^{\frac{18d}{9d-4}}_x} \\
& \qquad + \big\||\nabla|I u_H\big\|^2_{L^9_tL^{\frac{18d}{9d-4}}_x}
\big\||\nabla| Iu_L\big\|_{L^9_tL^{\frac{18d}{9d-4}}_x}\lesssim
\big\|Iu\big\|^3_{X^{1, \frac12+}}.
\endaligned
\end{equation*}

Now we estimate the fourth term. By Sobolev's inequality and
H\"{o}lder's inequality, we obtain
\begin{equation*}
\aligned \big\|\frac{1}{N_{1,2,3}} P_{N_{1,2,3}}\big( I  \big(u_H
|\nabla|^{-(d-\gamma)}&|u_L|^2\big)\big)\big\|_{L^3_tL^{\frac{6d}{3d-4}}_x}
\\
& \lesssim \big\| |\nabla|^{-1} P_{N_{1,2,3}}\big( I \big(u_H
|\nabla|^{-(d-\gamma)}|u_L|^2\big)\big)\big\|_{L^3_tL^{\frac{6d}{3d-4}}_x}
\\
& \lesssim \big\| u_H
|\nabla|^{-(d-\gamma)}|u_L|^2\big\|_{L^3_tL^{\frac{6d}{3d+2}}_x}
\\
& \lesssim \big\| u_H\big\|_{L^9_tL^{\frac{18d}{9d-4}}_x} \big\|
u_L\big\|^2_{L^9_tL^{\frac{18d}{9d-9\gamma+5}}_x} \lesssim
\big\|Iu\big\|^3_{X^{1, \frac12+}}.
\endaligned
\end{equation*}

The remainder terms are similar to the third and fourth terms
because we can ignore the complex conjugates. This completes the
proof.

Now we proceed to prove the sextilinear estimate. Note that in the
treatment of the quadrilinear form as in Proposition
\ref{increment1}, we always took the $\Delta u_1$ factor in
$L^3_tL^{6d/(3d-4)}_x$, estimating this by $N_1
\big\|Iu_1\big\|_{X^{1,\frac12+}}$. Together with Proposition
\ref{increment1} and Lemma \ref{trilinear}, we can obtain the
following estimate.

\begin{proposition}\label{increment2}
For any Schwartz function $u$, and any $\delta \thickapprox 1$ as in
Proposition \ref{MLWP}, we have that
\begin{equation*}
\aligned \Big| \int^{T+\delta}_{T}\Lambda_6\big(M_6; u \big)\ dt
\Big| \lesssim  N^{-1+} \big\|Iu \big\|^6_{X^{1, \frac12+}}
\endaligned
\end{equation*}
for $s>\frac{\gamma}{2}-1$.
\end{proposition}

\section{Almost Interaction Morawetz estimate } \setcounter{section}{4}\setcounter{equation}{0}
In this section, we aim to prove the interaction Morawetz estimate
for the smoothed out solution $Iu$, that is, ``almost Morawetz
estimate''. For this, we consider $a(x_1,x_2)=|x_1-x_2|:\mathbb{R}^d
\times \mathbb{R}^d \rightarrow \mathbb{R}$, a convex and locally
integrable function of polynomial growth. In all of our arguments,
we will work with the Schwarz solutions. This will simplify the
calculations and will enable us to justify the steps in the
subsequent proofs. Then we approximate the $H^s$ solutions by the
schwarz solutions.
\begin{theorem}\label{almostinteraction} Let
$u$ be a Schwarz solution to
\begin{equation*}
\aligned
    iu_t +  \Delta u  & = \widetilde{\mathcal{N}}(u),  \ (x, t) \in \mathbb{R}^d \times [0,T], \\
\endaligned
\end{equation*}
where $\widetilde{\mathcal{N}}(u)=\big( |x|^{-\gamma}* |u|^2 \big)
u$. Let $Iu $ be a solution to
\begin{equation}\label{iu}
\aligned
    iIu_t +  \Delta Iu  & = I\big(\widetilde{\mathcal{N}}(u)\big),  \ (x, t) \in \mathbb{R}^d \times
    [0,T].
\endaligned
\end{equation}
Then
\begin{equation}\label{AIME}
\aligned \big\| |\nabla|^{-\frac{d-3}{4}}I u\big\|^4_{L^4_TL^4_x}
\lesssim \big\|Iu&\big\|_{L^{\infty}_T\dot{H}^{1}_x}
\big\|Iu\big\|^3_{L^{\infty}_TL^2_x} \\
+&  \int^T_0 \int_{\mathbb{R}^d\times \mathbb{R}^d}\nabla a \cdot
\big\{\widetilde{\mathcal{N}}_{bad}, Iu(t, x_1)Iu(t,x_2) \big\}_{p}
dx_1 dx_2 dt.
\endaligned
\end{equation}
with $\big\{\cdot, \cdot \big\}_{p}$ is the momentum bracket defined
by
\begin{equation*}
\aligned \big\{f,g \big\}_{p}=\text{Re}\big(f\nabla \overline{g}
-g\nabla \overline{f}\big),
\endaligned
\end{equation*}
and
\begin{equation*} \aligned
\widetilde{\mathcal{N}}_{bad}=\sum^2_{i=1}\big(
I\widetilde{\mathcal{N}}_i(u_i)-\widetilde{\mathcal{N}}_i(Iu_i)
\big) \prod^2_{j=1, j\not=i}Iu_j,
\endaligned
\end{equation*}
where $u_i$ is a solution to
\begin{equation} \label{equ2}
  \aligned
    iu_t +  \Delta u  & = \widetilde{\mathcal{N}}(u),  \ (x_i, t) \in \mathbb{R}^d \times \mathbb{R}, \quad d\geq 3,\\
\endaligned
 \end{equation}
here $x_i \in \mathbb{R}^d$, not a coordinate. In particular, on a
time interval $J_k$ where the local well-posedness Proposition
\ref{MLWP} holds, we have that
\begin{equation*}
\aligned \int_{J_k} \int_{\mathbb{R}^d\times \mathbb{R}^d}\nabla a
\cdot \big\{\widetilde{\mathcal{N}}_{bad}, Iu(t, x_1)Iu(t,x_2)
\big\}_{p} dx_1 dx_2 dt \lesssim
\frac{1}{N^{1-}}\big\|u\big\|^6_{Z_I(J_k)}.
\endaligned
\end{equation*}
\end{theorem}

Toward this goal, we recall the idea of the proof of the interaction
Morawetz estimate for the defocusing nonlinear cubic Schr\"{o}dinger
equation in three space dimensions \cite{CKSTT04}. We present the
result using a tensor of Schr\"{o}dinger solutions that emerged in
\cite{CGT07}, \cite{DPST07}. We first recall the generalized Virial
identity
 \cite{CGT07}, \cite{LinS78}.
\begin{proposition}\label{im}
Let $a:\mathbb{R}^d \longrightarrow \mathbb{R}$ be convex and $u$ be
a smooth solution to the solution
\begin{equation}\label{u}
 \aligned
    iu_t +  \Delta u  & = \widetilde{\mathcal{N}}(u), \quad (t,x)\in [0,
    T]    \times \mathbb{R}^d.
\endaligned
\end{equation}
Then the following inequality holds
\begin{equation*}
\aligned \int^T_0 \int_{\mathbb{R}^d} \big( -\Delta \Delta
a\big)\big|u(t,x) \big|^2dxdt + 2\int^T_0 \int_{\mathbb{R}^d} \nabla
a \cdot \big\{ \widetilde{\mathcal{N}} , u \big\}_{p} dxdt \lesssim
\big|M_a(T)-M_a(0) \big|
\endaligned
\end{equation*}
where $M_a(t)$ is the Morawetz action corresponding to $u$ and is
given by
\begin{equation*}
\aligned M_a(t)=2\int_{\mathbb{R}^d}\nabla a(x) \cdot
\text{Im}\big(\overline{u}(x)\nabla u(x) \big) dx.
\endaligned
\end{equation*}
\end{proposition}

{\bf Proof of Theorem \ref{almostinteraction}: }Now we rewrite the
equation (\ref{iu}) as
\begin{equation*}
\aligned
 iIu_t +  \Delta Iu  & = \widetilde{\mathcal{N}}(Iu)+\big(
 I\big(\widetilde{\mathcal{N}}(u)\big)-\widetilde{\mathcal{N}}(Iu)\big),
\endaligned
\end{equation*}
then by symmetry, the term $\widetilde{\mathcal{N}}(Iu)$ will create
a positive term that we can ignore, which is the same to the case in
\cite{MiXZ07b}. While the commutator
$I\big(\widetilde{\mathcal{N}}(u)\big)-\widetilde{\mathcal{N}}(Iu)$
will introduce an error term.  Thus by Proposition \ref{im}, we have
\begin{equation*}
\aligned \int^T_0 \int \big( -\Delta \Delta a\big)\big|Iu(t,x)
\big|^2dxdt \lesssim  \sup_{t\in [0,T]} & \Big|\int \nabla a(x)
\cdot
\text{Im}\big(\overline{Iu}(x)\nabla Iu(x) \big) dx \Big| \\
+& \Big|\int^T_0 \int \nabla a \cdot \big\{
I\widetilde{\mathcal{N}}(u)-\widetilde{\mathcal{N}}(Iu) , Iu
\big\}_{p} dxdt\Big|.
\endaligned
\end{equation*}
The second term on the right hand side of this inequality is what we
call an error. We now turn to the details. The conjugation will play
no crucial role in the forthcoming argument.

Now define the tensor product $u:=\big(u_1\otimes u_2\big)(t,x)$ for
$x$ in
\begin{equation*}
\aligned \mathbb{R}^{d+d}=\big\{x=(x_1, x_2): x_1 \in \mathbb{R}^d,
x_2 \in \mathbb{R}^d \big\}
\endaligned
\end{equation*}
by the formula
\begin{equation*}
\aligned \big(u_1\otimes u_2\big)(t,x)=u_1(t,x_1)u_2(t,x_2).
\endaligned
\end{equation*}
let us set
\begin{equation*}
\aligned IU(t,x)=\prod^2_{j=1}Iu(t,x_j).
\endaligned
\end{equation*}
If $u$ solves (\ref{u}) for $d$ dimensions, then $IU$ solves
(\ref{u}) for $2d$ dimensions, with right hand side
$\widetilde{\mathcal{N}}_I$ given by
\begin{equation*}
\aligned \widetilde{\mathcal{N}}_I=\sum^2_{i=1}\Big(
I\big(\widetilde{\mathcal{N}}_i(u_i)\big) \prod^2_{j=1,
j\not=i}Iu_j\Big).
\endaligned
\end{equation*}
Now let us decompose
\begin{equation*}
\aligned
 \widetilde{\mathcal{N}}_I &=  \widetilde{\mathcal{N}}_{good} +
 \widetilde{\mathcal{N}}_{bad}\\
 &\triangleq \sum^2_{i=1}\Big(
\widetilde{\mathcal{N}}_i(Iu_i) \prod^2_{j=1, j\not=i}Iu_j\Big) +
\sum^2_{i=1}\Big( \big(I(\widetilde{\mathcal{N}}_i(u_i))-
\widetilde{\mathcal{N}}_i(Iu_i)\big) \prod^2_{j=1,
j\not=i}Iu_j\Big).
\endaligned
\end{equation*}
The first term summand creates a positive term that we can ignore
again. The term we call $ \widetilde{\mathcal{N}}_{bad}$ produces
the error term. Now we pick $a(x)=a(x_1, x_2)=|x_1-x_2|$ where
$(x_1,x_2)\in \mathbb{R}^d \times \mathbb{R}^d$. Hence we have
\begin{equation*}
\aligned \big\| |\nabla|^{-\frac{d-3}{4}}I u\big\|^4_{L^4_TL^4_x}
\lesssim \big\|Iu&\big\|_{L^{\infty}_T\dot{H}^{1}_x}
\big\|Iu\big\|^3_{L^{\infty}_TL^2_x} \\
+&  \Big|\int^T_0 \int_{\mathbb{R}^d\times \mathbb{R}^d}\nabla a
\cdot \big\{\widetilde{\mathcal{N}}_{bad}, Iu(t, x_1)Iu(t,x_2)
\big\}_{p} dx_1 dx_2 dt\Big|.
\endaligned
\end{equation*}
Note that the second term of the right hand side comes from the
momentum bracket term in the proof of Proposition \ref{im}.
Following with the same calculations in \cite{CGT07}, we deduce that
\begin{equation}\label{estimate}
\aligned \mathcal{E}:&=  \Big|\int^T_0 \int_{\mathbb{R}^d\times
\mathbb{R}^d}\nabla a \cdot \big\{\widetilde{\mathcal{N}}_{bad},
Iu(t, x_1)Iu(t,x_2) \big\}_{p} dx_1 dx_2 dt\Big|\\
&\lesssim \big(
\big\|I\big(\widetilde{\mathcal{N}}(u)\big)-\widetilde{\mathcal{N}}(Iu)\big\|_{L^1_tL^2_x}
+
\big\|\nabla_x\big(I\big(\widetilde{\mathcal{N}}(u)\big)-\widetilde{\mathcal{N}}(Iu)\big)\big\|_{L^1_tL^2_x}\big)
\big\|u\big\|^3_{Z_I(J)}.
\endaligned
\end{equation}

Now we proceed to estimate
$\big\|\nabla_x\big(I\big(\widetilde{\mathcal{N}}(u)\big)-\widetilde{\mathcal{N}}(Iu)\big)\big\|_{L^1_tL^2_x}$,
which is the harder term. The term
$\big\|I\big(\widetilde{\mathcal{N}}(u)\big)-\widetilde{\mathcal{N}}(Iu)\big\|_{L^1_tL^2_x}$
can be estimated in the same way. Note that
\begin{equation*}
\aligned \widetilde{\mathcal{N}}(u)=\big( |x|^{-\gamma}* |u|^2 \big)
u,
\endaligned
\end{equation*}
we have
\begin{equation*}
\aligned \mathcal{F}_x
\Big(\nabla_x\big(&I\big(\widetilde{\mathcal{N}}(u)\big)-\widetilde{\mathcal{N}}(Iu)\big)\Big)(\xi)\\
=&\int_{\xi=\sum^3_{j=1}\xi_j}
i\xi|\xi_{2,3}|^{-(d-\gamma)}\big(m(\xi)-m(\xi_1)m(\xi_2)m(\xi_3)\big)\widehat{u}(\xi_1)
\widehat{\overline{u}}(\xi_2)\widehat{u}(\xi_3)d\xi_1d\xi_2d\xi_3.
\endaligned
\end{equation*}

We decompose $u$ into a sum of dyadic pieces $u_j$ localized around
$N_j$, then
\begin{equation*}
\aligned &\  \big\|\nabla_x\big(
I\big(\widetilde{\mathcal{N}}(u)\big)-\widetilde{\mathcal{N}}(Iu)\big)\big\|_{L^1_tL^2_x}
\\
= & \ \big\|\mathcal{F}_x
\Big(\nabla_x\big(I\big(\widetilde{\mathcal{N}}(u)\big)-\widetilde{\mathcal{N}}(Iu)\big)\Big)(\xi)
\big\|_{L^1_tL^2_{\xi}} \\
\lesssim &\  \sum_{N_1, N_2, N_3} \Big\| \int_{|\xi_j|\thickapprox
N_j,\atop \xi=\sum^3_{j=1}\xi_j}
\big|\xi\big||\xi_{2,3}|^{-(d-\gamma)}\big|m(\xi)-m(\xi_1)m(\xi_2)m(\xi_3)\big|\\
& \qquad \qquad \qquad \qquad \qquad \qquad \qquad \qquad
\qquad\times \widehat{u}(\xi_1)
\widehat{\overline{u}}(\xi_2)\widehat{u}(\xi_3)d\xi_1d\xi_2d\xi_3\Big\|_{L^1_tL^2_{\xi}}.
\endaligned
\end{equation*}

Since the conjugation plays no crucial role here, without loss of
generality, we assume that
\begin{equation*}
\aligned N_1 \geq N_2 \geq N_3.
\endaligned
\end{equation*}

Set
\begin{equation*}
\aligned \sigma(\xi_1, \xi_2, \xi_3)=\big| \xi_1+\xi_2+ \xi_3\big|
\big|m(\xi_1+\xi_2+ \xi_3)-m(\xi_1)m(\xi_2)m(\xi_3)\big|,
\endaligned
\end{equation*}
then
\begin{equation*}
\aligned \sigma(\xi_1, \xi_2, \xi_3)&=\sum^4_{j=1}\chi_j(\xi_1,
\xi_2, \xi_3) \sigma(\xi_1, \xi_2, \xi_3)\\
:&=\sum^4_{j=1}\sigma_j(\xi_1, \xi_2, \xi_3),
\endaligned
\end{equation*}
where $\chi_j(\xi_1, \xi_2, \xi_3)$ is a smooth characteristic
function of the set $\Omega_j$ defined as follows:
\begin{enumerate}
\item[$\bullet$] $\Omega_1=\big\{ |\xi_i|\thickapprox N_i, i=1,2,3; N\gg N_1 \big\};$
\item[$\bullet$] $\Omega_2=\big\{ |\xi_i|\thickapprox N_i, i=1,2,3;  N_1\gtrsim N\gg N_2\big\};$
\item[$\bullet$] $\Omega_3=\big\{ |\xi_i|\thickapprox N_i, i=1,2,3;  N_1\geq N_2 \gtrsim N\gg N_3 \big\};$
\item[$\bullet$] $\Omega_4=\big\{ |\xi_i|\thickapprox N_i, i=1,2,3;  N_1\geq  N_2 \geq  N_3 \gtrsim N \big\}.$
\end{enumerate}
Hence, we have
\begin{equation*}
\aligned \big\|\nabla_x\big(
I\big(&\widetilde{\mathcal{N}}(u)\big)-\widetilde{\mathcal{N}}(Iu)\big)\big\|_{L^1_tL^2_x}
\\
 \lesssim &\sum_{N_1, N_2, N_3} \sum^4_{j=1} \Big\|
\int_{|\xi_j|\thickapprox N_j,\atop\xi=\sum^3_{j=1}\xi_j }
|\xi_{2,3}|^{-(d-\gamma)} \sigma_j(\xi_1, \xi_2,
\xi_3)\widehat{u}(\xi_1)
\widehat{\overline{u}}(\xi_2)\widehat{u}(\xi_3)d\xi_1d\xi_2d\xi_3\Big\|_{L^1_tL^2_{\xi}}
\\
:=& \sum_{N_1, N_2, N_3} \sum^4_{j=1} L_j.
\endaligned
\end{equation*}

{\bf Contribution of $L_1.$} Since $\sigma_1$ is identically zero
when $N \geq 4N_1$, $L_1$ gives no contribution to the sum above.

{\bf Contribution of $L_2.$} By the mean value theorem, we have the
pointwise bound
\begin{equation*}
\aligned \sigma_2(\xi_1, \xi_2, \xi_3) \lesssim N_1 \cdot m_1
\frac{N_2}{N_1} = m_1 N_2.
\endaligned
\end{equation*}
Hence, by H\"{o}lder's inequality and Hardy-Littlewood-Sobolev's
inequality, we obtain
\begin{equation*}
\aligned L_2 =&\Big\| \int_{|\xi_j|\thickapprox N_j,\atop\xi =
\sum^3_{j=1}\xi_j } |\xi_{2,3}|^{-(d-\gamma)} \sigma_2(\xi_1, \xi_2,
\xi_3)\widehat{u}(\xi_1)
\widehat{u}(\xi_2)\widehat{u}(\xi_3)d\xi_1d\xi_2d\xi_3\Big\|_{L^1_tL^2_{\xi}}\\
\lesssim & m_1 N_2 \Big\| \int_{|\xi_j|\thickapprox N_j,\atop\xi=
\sum^3_{j=1}\xi_j } |\xi_{2,3}|^{-(d-\gamma)} \widehat{u}(\xi_1)
\widehat{u}(\xi_2)\widehat{u}(\xi_3)d\xi_1d\xi_2d\xi_3\Big\|_{L^1_tL^2_{\xi}}\\
\lesssim & m_1 N_2
\big\|u_1\big\|_{L^{3}_tL^{\frac{6d}{3d-4}}_x}\big\|u_2\big\|_{L^{3}_tL^{\frac{6d}{3d-4}}_x}\big\|u_3\big\|_{L^{3}_tL^{\frac{6d}{3d-6\gamma+8}}_x}\\
\lesssim & m_1 N_2 N^{\gamma-2}_3 \prod^3_{J=1}
\big\|u_j\big\|_{X^{0,\frac12+}}.
\endaligned
\end{equation*}

It suffices to show that
\begin{equation*}
\aligned m_1 N_2 N^{\gamma-2}_3 \lesssim N^{-1+}N^{0-}_1
m_1N_1\langle N_2\rangle \langle N_3\rangle.
\endaligned
\end{equation*}
We reduce to show that
\begin{equation*}
\aligned N^{1-}N^{0+}_1   \lesssim  N_1\langle N_2\rangle N^{-1}_2
\langle N_3\rangle N^{2-\gamma}_3.
\endaligned
\end{equation*}
This is true since
\begin{equation*}
\aligned N_1  \gtrsim & N^{1-}N^{0+}_1 ;\\
\langle N_2\rangle N^{-1}_2  \gtrsim 1; &\quad  \langle N_3\rangle
N^{2-\gamma}_3 \gtrsim 1.
\endaligned
\end{equation*}

{\bf Contribution of $L_3.$} Note that
\begin{equation*}
\aligned \sigma_3(\xi_1, \xi_2, \xi_3) \lesssim N_1m_1 + N_1m_1m_2
\lesssim N_1m_1.
\endaligned
\end{equation*}
Hence, by H\"{o}lder's inequality and Hardy-Littlewood-Sobolev's
inequality, we have
\begin{equation*}
\aligned L_3 =&\Big\| \int_{|\xi_j|\thickapprox N_j,\atop\xi =
\sum^3_{j=1}\xi_j } |\xi_{2,3}|^{-(d-\gamma)} \sigma_3(\xi_1, \xi_2,
\xi_3)\widehat{u}(\xi_1)
\widehat{u}(\xi_2)\widehat{u}(\xi_3)d\xi_1d\xi_2d\xi_3\Big\|_{L^1_tL^2_{\xi}}\\
\lesssim & m_1 N_1\Big\| \int_{|\xi_j|\thickapprox N_j,\atop\xi=
\sum^3_{j=1}\xi_j } |\xi_{2,3}|^{-(d-\gamma)} \widehat{u}(\xi_1)
\widehat{u}(\xi_2)\widehat{u}(\xi_3)d\xi_1d\xi_2d\xi_3\Big\|_{L^1_tL^2_{\xi}}\\
\lesssim & m_1 N_1
\big\|u_1\big\|_{L^{3}_tL^{\frac{6d}{3d-4}}_x}\big\|u_2\big\|_{L^{3}_tL^{\frac{6d}{3d-4}}_x}\big\|u_3\big\|_{L^{3}_tL^{\frac{6d}{3d-6\gamma+8}}_x}\\
\lesssim & m_1 N_1 N^{\gamma-2}_3 \prod^3_{J=1}
\big\|u_j\big\|_{X^{0,\frac12+}}.
\endaligned
\end{equation*}

It suffices to show that
\begin{equation*}
\aligned m_1 N_1 N^{\gamma-2}_3 \lesssim N^{-1+}N^{0-}_2 m_1N_1m_2
N_2 \langle N_3\rangle.
\endaligned
\end{equation*}
We reduce to show that
\begin{equation*}
\aligned N^{1-}N^{0+}_2   \lesssim  m_2 N_2  \langle  N_3\rangle
N^{2-\gamma}_3.
\endaligned
\end{equation*}
This is true since
\begin{equation*}
\aligned m_2N_2  \gtrsim  N^{1-}N^{0+}_2 ;\quad
 \langle N_3\rangle
N^{2-\gamma}_3 \gtrsim  1.
\endaligned
\end{equation*}

{\bf Contribution of $L_4.$} Note that
\begin{equation*}
\aligned \sigma_4(\xi_1, \xi_2, \xi_3) \lesssim N_1m_1 + N_1m_1m_2
\lesssim N_1m_1.
\endaligned
\end{equation*}
Hence, by H\"{o}lder's inequality and Hardy-Littlewood-Sobolev's
inequality, we obtain
\begin{equation*}
\aligned L_4 =&\Big\| \int_{|\xi_j|\thickapprox N_j,\atop\xi =
\sum^3_{j=1}\xi_j } |\xi_{2,3}|^{-(d-\gamma)} \sigma_4(\xi_1, \xi_2,
\xi_3)\widehat{u}(\xi_1)
\widehat{u}(\xi_2)\widehat{u}(\xi_3)d\xi_1d\xi_2d\xi_3\Big\|_{L^1_tL^2_{\xi}}\\
\lesssim & m_1 N_1\Big\| \int_{|\xi_j|\thickapprox N_j,\atop\xi=
\sum^3_{j=1}\xi_j } |\xi_{2,3}|^{-(d-\gamma)} \widehat{u}(\xi_1)
\widehat{u}(\xi_2)\widehat{u}(\xi_3)d\xi_1d\xi_2d\xi_3\Big\|_{L^1_tL^2_{\xi}}\\
\lesssim & m_1 N_1
\big\|u_1\big\|_{L^{3}_tL^{\frac{6d}{3d-4}}_x}\big\|u_2\big\|_{L^{3}_tL^{\frac{6d}{3d-4}}_x}\big\|u_3\big\|_{L^{3}_tL^{\frac{6d}{3d-6\gamma+8}}_x}\\
\lesssim & m_1 N_1 N^{\gamma-2}_3 \prod^3_{J=1}
\big\|u_j\big\|_{X^{0,\frac12+}}.
\endaligned
\end{equation*}

It suffices to show that
\begin{equation*}
\aligned m_1 N_1 N^{\gamma-2}_3 \lesssim N^{-1+}N^{0-}_2 m_1N_1m_2
N_2 m_3 N_3.
\endaligned
\end{equation*}
We reduce to show that
\begin{equation*}
\aligned N^{1-}N^{0+}_2   \lesssim  m_2 N_2  m_3  N_3
N^{2-\gamma}_3.
\endaligned
\end{equation*}
This is true since for $s\geq \gamma-2$, we have
\begin{equation*}
\aligned m_2  N_2\ m_3 N^{3-\gamma}_3 & \gtrsim m_2 N_2 \gtrsim
N^{1-}N^{0+}_2
\endaligned
\end{equation*}
where we used the fact that $m(\xi)\langle\xi\rangle^p$ is monotone
non-decreasing if $s+p\geq 1$. While for $\frac{\gamma}{2}-1 < s <
\gamma-2$, we have
\begin{equation*}
\aligned m_2  N_2\ m_3 N^{3-\gamma}_3 & \gtrsim
m_2  N_2\ m_2 N^{3-\gamma}_2\\
& \gtrsim N^{4-\gamma-}N^{0+}_2 \gtrsim N^{1-}N^{0+}_2
\endaligned
\end{equation*}
where we used the fact that $m(\xi)\langle\xi\rangle^p$ is monotone
non-increasing if $s+p < 1$.

\section{Proof of Theorem \ref{theorem}}
\setcounter{section}{5} \setcounter{equation}{0} We first scale the
solution. Suppose that $u(t,x)$ is a global in time solution to
(\ref{equ1}) with initial data $u_0 \in C^{\infty}_0\big(
\mathbb{R}^d\big).$ Setting
\begin{equation*}
\aligned
u^{\lambda}(t,x)=\lambda^{-\frac{n+2-\gamma}{2}}u(\frac{t}{\lambda^2},
\frac{x}{\lambda}),
\endaligned
\end{equation*}
we choose a parameter $\lambda$ so that
$\big\|Iu^{\lambda}_0\big\|_{H^1}=O(1)$, that is
\begin{equation}\label{namda}
\aligned \lambda \approx N^{\frac{1-s}{s-\frac{\gamma}{2}+1}}.
\endaligned
\end{equation}

Next, let us define
\begin{equation*}
\aligned S:=\big\{0\leq t<\infty:
\big\|Iu^{\lambda}\big\|_{L^4\dot{H}^{-\frac{d-3}{4},4}\big(
[0,t]\times \mathbb{R}^d\big)}\leq K
\lambda^{\frac{3}{4}(\frac{\gamma}{2}-1)} \big\},
\endaligned
\end{equation*}
with $K$ a constant to be chosen later. We claim that $S$ is the
whole interval $[0, \infty)$.  Indeed, assume by contradiction that
it is not so, then since
\begin{equation*}
\aligned
\big\|Iu^{\lambda}\big\|_{L^4\dot{H}^{-\frac{d-3}{4},4}\big(
[0,t]\times \mathbb{R}^d\big)}
\endaligned
\end{equation*}
is a continuous function of time, there exists a time $T\in [0,
\infty)$ such that
\begin{eqnarray}
\big\|Iu^{\lambda}\big\|_{L^4\dot{H}^{-\frac{d-3}{4},4}\big(
[0,T]\times \mathbb{R}^d\big)} & > & K
\lambda^{\frac{3}{4}(\frac{\gamma}{2}-1)}, \label{contradiciton}\\
\big\|Iu^{\lambda}\big\|_{L^4\dot{H}^{-\frac{d-3}{4},4}\big(
[0,T]\times \mathbb{R}^d\big)} & \leq & 2 K
\lambda^{\frac{3}{4}(\frac{\gamma}{2}-1)}. \label{bound3}
\end{eqnarray}

We now split the interval $[0, T]$ into subintervals $J_k,
k=1,\cdots, L$ in such a way that
\begin{equation*}
\aligned
\big\|Iu^{\lambda}\big\|^4_{L^4\dot{H}^{-\frac{d-3}{4},4}\big(
J_k\times \mathbb{R}^d\big)}\leq \mu,
\endaligned
\end{equation*}
with $\mu$ as in Proposition \ref{universalbound}. This is possible
because of (\ref{bound3}). Then, the number $L$ of possible
subintervals must satisfy
\begin{equation}\label{iteratenumber}
\aligned L\approx \frac{\big(2 K
\lambda^{\frac{3}{4}(\frac{\gamma}{2}-1)}\big)^4}{\mu} \approx
\frac{\big(2 K\big)^4 \lambda^{ 3(\frac{\gamma}{2}-1)}}{\mu}.
\endaligned
\end{equation}
From Proposition \ref{MLWP} and Proposition \ref{universalbound}, we
know that
\begin{equation*}
\aligned \sup_{t\in[0,T]}E\big(Iu^{\lambda}(t)\big) \lesssim
E\big(Iu^{\lambda}_0\big) + \frac{L}{N^{1-}}
\endaligned
\end{equation*}
and by our choice (\ref{namda}) of $\lambda$,
$E\big(Iu^{\lambda}_0\big) \lesssim 1$. Hence, in order to guarantee
that
\begin{equation*}
\aligned E\big(Iu^{\lambda}(t)\big) \lesssim 1
\endaligned
\end{equation*}
holds for all $t\in [0,T]$, we need to require that
\begin{equation*}
\aligned L \lesssim N^{1-}.
\endaligned
\end{equation*}
According to (\ref{iteratenumber}), this is fulfilled as long as
\begin{equation}\label{con}
\aligned \frac{\big(2 K\big)^4 \lambda^{
3(\frac{\gamma}{2}-1)}}{\mu} \lesssim N^{1-}.
\endaligned
\end{equation}
From our choice of $\lambda$, the expression (\ref{con}) implies
that
\begin{equation*}
\aligned \frac{\big(2 K\big)^4 }{\mu} \lesssim N^{1-
\frac{1-s}{s-\frac{\gamma}{2}+1}3(\frac{\gamma}{2}-1)- }.
\endaligned
\end{equation*}
Thus this is possible for $s> \frac{4(\gamma-2)}{3\gamma-4}$ and a
large number $N$.

Now recall the a priori estimate (\ref{AIME})
\begin{equation*}
\aligned \big\| |\nabla|^{-\frac{d-3}{4}}I
u^{\lambda}\big\|^4_{L^4_TL^4_x} \lesssim
\big\|Iu^{\lambda}&\big\|_{L^{\infty}_T\dot{H}^{1}_x}
 \big\|Iu^{\lambda}\big\|^3_{L^{\infty}_TL^2_x} \\
+& \int^T_0 \int_{\mathbb{R}^d\times \mathbb{R}^d}\nabla a \cdot
\big\{\widetilde{\mathcal{N}}_{bad}, Iu^{\lambda}(t,
x_1)Iu^{\lambda}(t,x_2) \big\}_{p} dx_1 dx_2 dt.
\endaligned
\end{equation*}
Set
\begin{equation*}
\aligned Error(t):=\int_{\mathbb{R}^d\times \mathbb{R}^d}\nabla a
\cdot \big\{\widetilde{\mathcal{N}}_{bad}, Iu^{\lambda}(t,
x_1)Iu^{\lambda}(t,x_2) \big\}_{p} dx_1 dx_2.
\endaligned
\end{equation*}
By Theorem \ref{almostinteraction} and Proposition
\ref{universalbound} on each interval $J_k$, we have that
\begin{equation*}
\aligned \Big|\int_{J_k} Error(t)dt \Big| \lesssim \frac{1}{N^{1-}}
\big\|u^{\lambda}\big\|^6_{Z_I(J_k)} \lesssim \frac{1}{N^{1-}}.
\endaligned
\end{equation*}
Summing all the $J_k$'s, we have that
\begin{equation*}
\aligned \Big|\int^T_{0} Error(t)dt \Big|  \lesssim \frac{L}{N^{1-}}
\lesssim N^{0+}.
\endaligned
\end{equation*}
Therefore, by our choice (\ref{namda}) of $\lambda$, we obtain
\begin{equation*}
\aligned \big\| |\nabla|^{-\frac{d-3}{4}}
Iu^{\lambda}\big\|^4_{L^4_TL^4_x}&  \lesssim
\big\|Iu^{\lambda}\big\|_{L^{\infty}_T\dot{H}^{1}_x}
\big\|Iu^{\lambda}\big\|^3_{L^{\infty}_TL^2_x} +N^{0+} \\
& \lesssim C \lambda^{3(\frac{\gamma}{2}-1)}.
\endaligned
\end{equation*}
This estimate contradicts (\ref{contradiciton}) for an appropriate
choice of $K$. Hence $S=[0, \infty)$. In addition, let $T_0$  be
chosen arbitrarily, we have also proved that for $s>
\frac{4(\gamma-2)}{3\gamma-4}$,
\begin{equation*}
\aligned \big\|Iu^{\lambda}(\lambda^2T_0)\big\|_{H^1_x}=O(1).
\endaligned
\end{equation*}
Then
\begin{equation*}
\aligned \big\|u(T_0)\big\|_{H^s}& = \big\|u(T_0)\big\|_{L^2} +
\big\|u(T_0)\big\|_{\dot{H}^s}\\
& = \big\|u_0\big\|_{L^2} +
\lambda^{s-\frac{\gamma}{2}+1} \big\|u^{\lambda}(\lambda^2T_0)\big\|_{\dot{H}^s}\\
& \lesssim
\lambda^{s-\frac{\gamma}{2}+1} \big\|I u^{\lambda}(\lambda^2T_0)\big\|_{H^1} \lesssim \lambda^{s-\frac{\gamma}{2}+1}\approx N^{1-s}.\\
\endaligned
\end{equation*}
Since $T_0$ is arbitrarily large, the a priori bound on the $H^s$
norm concludes the global well-posednes of the Cauchy problem
(\ref{equ1}).

Note that we have obtained that
\begin{equation*}
\aligned \big\|Iu\big\|_{L^4\dot{H}^{-\frac{d-3}{4},4}\big(
[0,+\infty)\times \mathbb{R}^d\big)}\leq C(\big\|u_0\big\|_{H^s}),
\endaligned
\end{equation*}
this together with Proposition \ref{linear}, Proposition
\ref{universalbound} and the property of the operator $I$ implies
that
\begin{equation*}
\aligned \sup_{(q,r)\; admissible} \big\|\langle \nabla \rangle^s
u\big\|_{L^qL^r\big([0,+\infty)\times \mathbb{R}^d\big)} \leq
\big\|u\big\|_{Z_I([0,+\infty))}\lesssim C\big(
\big\|Iu_0\big\|_{H^1}\big)\lesssim
C\big(\big\|u_0\big\|_{H^s}\big),
\endaligned
\end{equation*}
then we can prove scattering by using the well-known standard
argument \cite{Ca03}, \cite{CKSTT04} etc.. This completes the
proof.

\textbf{Acknowledgements:}   C. Miao  and G.Xu  were partly
supported by the NSF of China (No.10725102, No.10726053), and
L.Zhao was supported by China postdoctoral science foundation
project. The second author would like to thank Nikolaos Tzirakis
for helpful communications.

\begin{center}

\end{center}
\end{document}